\documentclass{amsart}
\usepackage{amsfonts}

\usepackage{graphicx}
\usepackage{amscd}

\newtheorem{theorem}{Theorem}
\theoremstyle{plain}

\newtheorem{corollary}{Corollary}

\newtheorem{lemma}{Lemma}

\newtheorem{remark}{Remark}

\numberwithin{equation}{section}

\input{tcilatex}

\begin{document}
\title[\v{C}eby\v{s}ev Functional]{Bounding the \v{C}eby\v{s}ev Functional for Sequences of Vectors in Normed
Linear Spaces}
\author{S.S. Dragomir}
\address{School of Computer Science \& Mathematics\\
Victoria University, Melbourne, Australia}
\email{sever@matilda.vu.edu.au}
\urladdr{http://rgmia.vu.edu.au/SSDragomirWeb.html}
\date{March 13, 2003}

\begin{abstract}
Some new bounds for \v{C}eby\v{s}ev functional for sequences of vectors in
normed linear spaces are pointed out.
\end{abstract}

\maketitle

\section{Introduction}

Consider the \v{C}eby\v{s}ev functional defined for $\mathbf{p}=\left(
p_{1},...,p_{n}\right) \in \mathbb{R}^{n},$ $\mathbf{\alpha }=\left( \alpha
_{1},...,\alpha _{n}\right) \in \mathbb{K}^{n}\left( \mathbb{K=R}\text{ or }%
\mathbb{C}\right) $ and $\mathbf{x}=\left( x_{1},...,x_{n}\right) \in X^{n},$
where $X$ is a linear space over the real or complex number field $\mathbb{K}
$:
\begin{equation}
T_{n}\left( \mathbf{p};\mathbf{\alpha },\mathbf{x}\right)
:=P_{n}\sum_{i=1}^{n}p_{i}\alpha _{i}x_{i}-\sum_{i=1}^{n}p_{i}\alpha
_{i}\cdot \sum_{i=1}^{n}p_{i}x_{i},  \label{1.1}
\end{equation}
where $P_{n}:=\sum_{i=1}^{n}p_{i}.$

The following Gr\"{u}ss type inequalities for sequences in normed linear
spaces hold.

\begin{theorem}
\label{t1.1}Let $\left( X,\left\Vert .\right\Vert \right) $ be a normed
linear space over the real or complex number field $\mathbb{K}$, $\mathbf{%
\alpha }=\left( \alpha _{1},...,\alpha _{n}\right) \in \mathbb{K}^{n},%
\mathbf{p}=\left( p_{1},...,p_{n}\right) \in \mathbb{R}_{+}^{n}$ with $%
\sum_{i=1}^{n}p_{i}=1$ and $\mathbf{x}=\left( x_{1},...,x_{n}\right) \in
X^{n}.$ Then one has the inequalities
\begin{eqnarray}
&&\left\Vert T_{n}\left( \mathbf{p};\mathbf{\alpha },\mathbf{x}\right)
\right\Vert  \label{1.2} \\
&\leq &\left\{ 
\begin{array}{l}
\left[ \sum_{i=1}^{n}i^{2}p_{i}-\left( \sum_{i=1}^{n}ip_{i}\right) ^{2}%
\right] \max_{1\leq j\leq n-1}\left\vert \Delta \alpha _{j}\right\vert
\max_{1\leq j\leq n-1}\left\Vert \Delta x_{j}\right\Vert ,\text{\cite{DB};}
\\ 
\\ 
\frac{1}{2}\sum_{i=1}^{n}p_{i}\left( 1-p_{i}\right)
\sum_{j=1}^{n-1}\left\vert \Delta \alpha _{j}\right\vert
\sum_{j=1}^{n-1}\left\Vert \Delta x_{j}\right\Vert ,\text{\cite{SSD2};} \\ 
\\ 
\sum_{1\leq i<j\leq n}p_{i}p_{j}\left( j-i\right) \left(
\sum_{j=1}^{n-1}\left\vert \Delta \alpha _{j}\right\vert ^{p}\right)
^{1/p}\left( \sum_{j=1}^{n-1}\left\Vert \Delta x_{j}\right\Vert ^{q}\right)
^{1/q}, \\ 
\\ 
p>1,\frac{1}{p}+\frac{1}{q}=1,\text{\cite{SSD1}}.
\end{array}
\right.  \notag
\end{eqnarray}
The constant $1$ in the first branch, $\frac{1}{2}$ in the second branch and 
$1$ in the third branch are best possible in the sense that they cannot be
replaced by smaller constants.
\end{theorem}

The following particular inequalities for unweighted means hold as well,
where $T_{n}\left( \mathbf{\alpha },\mathbf{x}\right) $ is defined as
follows:
\begin{equation*}
T_{n}\left( \mathbf{\alpha },\mathbf{x}\right) :=\frac{1}{n}%
\sum_{i=1}^{n}\alpha _{i}x_{i}-\frac{1}{n}\sum_{i=1}^{n}a_{i}\cdot \frac{1}{n%
}\sum_{i=1}^{n}x_{i}.
\end{equation*}

\begin{corollary}
\label{c1.2} With the assumptions of Theorem \ref{t1.1} for $X,\mathbf{%
\alpha }$ and $\mathbf{x}$, we have
\begin{eqnarray}
&&\left\Vert T_{n}\left( \mathbf{\alpha },\mathbf{x}\right) \right\Vert
\label{1.3} \\
&\leq &\left\{ 
\begin{array}{l}
\frac{1}{12}\left( n^{2}-1\right) \max_{1\leq j\leq n-1}\left\vert \Delta
\alpha _{j}\right\vert \max_{1\leq j\leq n-1}\left\Vert \Delta
x_{j}\right\Vert ,\text{\cite{DB};} \\ 
\\ 
\frac{1}{2}\cdot \left( 1-\frac{1}{n}\right) \sum_{j=1}^{n-1}\left\vert
\Delta \alpha _{j}\right\vert \sum_{j=1}^{n-1}\left\Vert \Delta
x_{j}\right\Vert ,\text{\cite{SSD2};} \\ 
\\ 
\frac{1}{6}\frac{n^{2}-1}{n}\left( \sum_{j=1}^{n-1}\left\vert \Delta \alpha
_{j}\right\vert ^{p}\right) ^{1/p}\left( \sum_{j=1}^{n-1}\left\Vert \Delta
x_{j}\right\Vert ^{q}\right) ^{1/q}, \\ 
\\ 
p>1,\frac{1}{p}+\frac{1}{q}=1,\text{\cite{SSD1}}.
\end{array}
\right.  \notag
\end{eqnarray}
Here the constants $\frac{1}{12},\frac{1}{2}$ and $\frac{1}{6}$ are best
possible in the sense that they cannot be replaced by smaller constants.
\end{corollary}

For applications to estimate the $p$-moments of guessing mappings, see \cite
{DB}. For applications in approximating the discrete Fourier transform, the
discrete Mellin transform as well as some applications for polynomials and
Lipschitzian mappings, see \cite{SSD1} and \cite{SSD2}.

For classical results related the \v{C}eby\v{s}ev functional, see \cite{CE1}%
, \cite{CE2}, \cite{HLP}, \cite{BIE}, \cite{SO}, \cite{KOR} and \cite{MP}.
For more recent results, see \cite{MP}, \cite{5b}, \cite{6b}, \cite{7b}, 
\cite{MPF1} and \cite{SSD}.

\section{The Identities}

The first result is embodied in the following

\begin{theorem}
\label{t2.1} Let $\mathbf{p}=\left( p_{1},...,p_{n}\right) ,\mathbf{a}%
=\left( a_{1},...,a_{n}\right) $ be $n$-tuples of real or complex numbers
and $\mathbf{x}=\left( x_{1},...,x_{n}\right) $ an $n$-tuple of vectors in
the linear space $X$. If we define 
\begin{eqnarray*}
P_{i} &:&=\sum_{k=1}^{i}p_{k},\bar{P}_{i}:=P_{n}-P_{i},i\in \left\{
1,...,n-1\right\} , \\
A_{i}\left( \mathbf{p}\right)  &:&=\sum_{k=1}^{i}p_{k}a_{k},\bar{A}%
_{i}\left( \mathbf{p}\right) :=A_{n}\left( \mathbf{p}\right) -A_{i}\left( 
\mathbf{p}\right) ,i\in \left\{ 1,...,n-1\right\} ,
\end{eqnarray*}
then we have the identity 
\begin{eqnarray}
&&T_{n}\left( \mathbf{p};\mathbf{a},\mathbf{x}\right) =\sum_{i=1}^{n-1}\det
\left( 
\begin{array}{cc}
P_{i} & P_{n} \\ 
A_{i}\left( \mathbf{p}\right)  & A_{n}\left( \mathbf{p}\right) 
\end{array}
\right) \cdot \Delta x_{i}  \label{2.1.a} \\
&=&P_{n}\sum_{i=1}^{n-1}P_{i}\left( \frac{A_{n}\left( \mathbf{p}\right) }{%
P_{n}}-\frac{A_{i}\left( \mathbf{p}\right) }{P_{i}}\right) \cdot \Delta
x_{i}\left( \text{if }P_{i}\neq 0,i\in \left\{ 1,...,n\right\} \right)  
\notag \\
&=&\sum_{i=1}^{n-1}P_{i}\bar{P}_{i}\left( \frac{\bar{A}_{i}\left( \mathbf{p}%
\right) }{\bar{P}_{i}}-\frac{A_{i}\left( \mathbf{p}\right) }{P_{i}}\right)
\cdot \Delta x_{i}\left( \text{if }P_{i},\bar{P}_{i}\neq 0,i\in \left\{
1,...,n-1\right\} \right) ;  \notag
\end{eqnarray}
where $\Delta x_{i}:=x_{i+1}-x_{i}\left( i\in \left\{ 1,...,n-1\right\}
\right) $ is the forward difference.
\end{theorem}

\begin{proof}
We use the following well known summation by parts formula 
\begin{equation}
\sum_{l=p}^{q-1}d_{l}\Delta v_{l}=\left. d_{l}v_{l}\right\vert
_{p}^{q}-\sum_{l=p}^{q-1}v_{l+1}\Delta d_{l},  \label{2.2}
\end{equation}
where $d_{l}$ are real or complex numbers, and $v_{l}$ are vectors in a
linear space, $l=p,...,q$ $\left( q>p;p,q\text{ are natural numbers}\right)
. $

If we choose in $\left( \ref{2.2}\right) $, $p=1,q=n,d_{i}=P_{i}A_{n}\left( 
\mathbf{p}\right) -P_{n}A_{i}\left( \mathbf{p}\right) $ and $v_{i}=x_{i}$ $%
\left( i\in \left\{ 1,...,n-1\right\} \right) ,$ then we get 
\begin{eqnarray*}
&&\sum_{i=1}^{n-1}\left( P_{i}A_{n}\left( \mathbf{p}\right)
-P_{n}A_{i}\left( \mathbf{p}\right) \right) \cdot \Delta x_{i} \\
&=&\left. \left[ P_{i}A_{n}\left( \mathbf{p}\right) -P_{n}A_{i}\left( 
\mathbf{p}\right) \right] \cdot x_{i}\right| _{1}^{n}-\sum_{i=1}^{n-1}\Delta
\left( P_{i}A_{n}\left( \mathbf{p}\right) -P_{n}A_{i}\left( \mathbf{p}%
\right) \right) \cdot x_{i+1} \\
&=&\left[ P_{n}A_{n}\left( \mathbf{p}\right) -P_{n}A_{n}\left( \mathbf{p}%
\right) \right] \cdot x_{n}-\left[ P_{1}A_{n}\left( \mathbf{p}\right)
-P_{n}A_{1}\left( \mathbf{p}\right) \right] \cdot x_{1} \\
&&-\sum_{i=1}^{n-1}\left[ P_{i+1}A_{n}\left( \mathbf{p}\right)
-P_{n}A_{i+1}\left( \mathbf{p}\right) -P_{i}A_{n}\left( \mathbf{p}\right)
+P_{n}A_{i}\left( \mathbf{p}\right) \right] \cdot x_{i+1} \\
&=&P_{n}p_{1}a_{1}x_{1}-p_{1}A_{n}\left( \mathbf{p}\right)
x_{1}-\sum_{i=1}^{n-1}\left( p_{i+1}A_{n}\left( \mathbf{p}\right)
-P_{n}p_{i+1}a_{i+1}\right) \cdot x_{i+1} \\
&=&P_{n}p_{1}a_{1}x_{1}-p_{1}A_{n}\left( \mathbf{p}\right) x_{1}-A_{n}\left( 
\mathbf{p}\right)
\sum_{i=1}^{n-1}p_{i+1}x_{i+1}+P_{n}\sum_{i=1}^{n-1}p_{i+1}a_{i+1}x_{i+1} \\
&=&P_{n}\sum_{i=1}^{n}p_{i}a_{i}x_{i}-\sum_{i=1}^{n}p_{i}a_{i}\cdot
\sum_{i=1}^{n}p_{i}x_{i} \\
&=&T_{n}\left( \mathbf{p};\mathbf{a},\mathbf{x}\right) ,
\end{eqnarray*}
which produce the first identity in $\left( \ref{2.1.a}\right) .$

The second and the third identities are obvious and we omit the details.
\end{proof}

Before we prove the second result, we need the following lemma providing an
identity that is interesting in itself as well.

\begin{lemma}
\label{l1.2} Let $\mathbf{p}=\left( p_{1},...,p_{n}\right) $ and $\mathbf{a}%
=\left( a_{1},...,a_{n}\right) $ be $n$-tuples of real or complex numbers.
Then we have the equality 
\begin{equation}
\det \left( 
\begin{array}{cc}
P_{i} & P_{n} \\ 
A_{i}\left( \mathbf{p}\right) & A_{n}\left( \mathbf{p}\right)
\end{array}
\right) =\sum_{j=1}^{n-1}P_{\min \left\{ i,j\right\} }\bar{P}_{\max \left\{
i,j\right\} }\cdot \Delta a_{j},  \label{2.3}
\end{equation}
for each $i\in \left\{ 1,...,n-1\right\} .$
\end{lemma}

\begin{proof}
Define, for $i\in \left\{ 1,...,n-1\right\} ,$ 
\begin{equation*}
K\left( i\right) :=\sum_{j=1}^{n-1}P_{\min \left\{ i,j\right\} }\bar{P}%
_{\max \left\{ i,j\right\} }\cdot \Delta a_{j}.
\end{equation*}
We have
\begin{eqnarray}
K\left( i\right)  &=&\sum_{j=1}^{i}P_{\min \left\{ i,j\right\} }\bar{P}%
_{\max \left\{ i,j\right\} }\cdot \Delta a_{j}+\sum_{j=i+1}^{n-1}P_{\min
\left\{ i,j\right\} }\bar{P}_{\max \left\{ i,j\right\} }\cdot \Delta a_{j}
\label{2.4.a} \\
&=&\sum_{j=1}^{i}P_{j}\bar{P}_{i}\cdot \Delta a_{j}+\sum_{j=i+1}^{n-1}P_{i}%
\bar{P}_{j}\cdot \Delta a_{j}  \notag \\
&=&\bar{P}_{i}\sum_{j=1}^{i}P_{j}\cdot \Delta a_{j}+P_{i}\sum_{j=i+1}^{n-1}%
\bar{P}_{j}\cdot \Delta a_{j}.  \notag
\end{eqnarray}
Using the summation by parts formula, we have
\begin{eqnarray}
\sum_{j=1}^{i}P_{j}\cdot \Delta a_{j} &=&\left. P_{j}\cdot a_{j}\right|
_{1}^{i+1}-\sum_{j=1}^{i}\left( P_{j+1}-P_{j}\right) \cdot a_{j+1}
\label{2.5} \\
&=&P_{i+1}a_{i+1}-p_{1}a_{1}-\sum_{j=1}^{i}p_{j+1}\cdot a_{j+1}  \notag \\
&=&P_{i+1}a_{i+1}-\sum_{j=1}^{i+1}p_{j}\cdot a_{j}  \notag
\end{eqnarray}
and 
\begin{eqnarray}
\sum_{j=i+1}^{n-1}\bar{P}_{j}\cdot \Delta a_{j} &=&\left. \bar{P}_{j}\cdot
a_{j}\right| _{i+1}^{n}-\sum_{j=i+1}^{n-1}\left( \bar{P}_{j+1}-\bar{P}%
_{j}\right) \cdot a_{j+1}  \label{2.6} \\
&=&\bar{P}_{n}a_{n}-\bar{P}_{i+1}a_{i+1}-\sum_{j=i+1}^{n-1}\left(
P_{n}-P_{j+1}-P_{n}+P_{j}\right) \cdot a_{j+1}  \notag \\
&=&-\bar{P}_{i+1}a_{i+1}+\sum_{j=i+1}^{n-1}p_{j+1}\cdot a_{j+1}.  \notag
\end{eqnarray}
Using $\left( \ref{2.5}\right) $ and $\left( \ref{2.6}\right) $ we have 
\begin{eqnarray*}
K\left( i\right)  &=&\bar{P}_{i}\left(
P_{i+1}a_{i+1}-\sum_{j=1}^{i+1}p_{j}\cdot a_{j}\right) +P_{i}\left(
\sum_{j=i+1}^{n-1}p_{j+1}\cdot a_{j+1}-\bar{P}_{i+1}a_{i+1}\right)  \\
&=&\bar{P}_{i}P_{i+1}a_{i+1}-P_{i}\bar{P}_{i+1}a_{i+1}-\bar{P}%
_{i}\sum_{j=1}^{i+1}p_{j}\cdot a_{j}+P_{i}\sum_{j=i+1}^{n-1}p_{j+1}\cdot
a_{j+1} \\
&=&\left[ \left( P_{n}-P_{i}\right) P_{i+1}-P_{i}\left( P_{n}-P_{i+1}\right) %
\right] a_{i+1} \\
&&+P_{i}\sum_{j=i+1}^{n-1}p_{j+1}\cdot a_{j+1}-\bar{P}_{i}%
\sum_{j=1}^{i+1}p_{j}\cdot a_{j}
\end{eqnarray*}
\begin{eqnarray*}
&=&P_{n}p_{i+1}a_{i+1}+P_{i}\sum_{j=i+1}^{n-1}p_{j+1}\cdot a_{j+1}-\bar{P}%
_{i}\sum_{j=1}^{i+1}p_{j}\cdot a_{j} \\
&=&\left( P_{i}+\bar{P}_{i}\right)
p_{i+1}a_{i+1}+P_{i}\sum_{j=i+1}^{n-1}p_{j+1}\cdot a_{j+1}-\bar{P}%
_{i}\sum_{j=1}^{i+1}p_{j}\cdot a_{j} \\
&=&P_{i}\sum_{j=i+1}^{n-1}p_{j}\cdot a_{j}-\bar{P}_{i}\sum_{j=1}^{i}p_{j}%
\cdot a_{j}=P_{i}\bar{A}_{i}\left( \mathbf{p}\right) -\bar{P}_{i}A_{i}\left( 
\mathbf{p}\right)  \\
&=&\det \left( 
\begin{array}{cc}
P_{i} & P_{n} \\ 
A_{i}\left( \mathbf{p}\right)  & A_{n}\left( \mathbf{p}\right) 
\end{array}
\right) ;
\end{eqnarray*}
and the identity is proved.
\end{proof}

We are able now to state and prove the second identity for the \v{C}eby\v{s}%
ev functional

\begin{theorem}
\label{t2.2} With the assumptions of Theorem \ref{t2.1}, we have the
equality 
\begin{equation}
T_{n}\left( \mathbf{p};\mathbf{a},\mathbf{x}\right)
=\sum_{i=1}^{n-1}\sum_{j=1}^{n-1}P_{\min \left\{ i,j\right\} }\bar{P}_{\max
\left\{ i,j\right\} }\cdot \Delta a_{j}\cdot \Delta x_{i}.  \label{2.7}
\end{equation}
\end{theorem}

The proof is obvious by Theorem \ref{t2.1} and Lemma \ref{l1.2}.

\begin{remark}
The identity $\left( \ref{2.7}\right) $, for n-tuples of real numbers, was
stated without a proof in paper \cite{MP}. It also may be found for the same
sequences in \cite[p. 281]{MPF1}, again without a proof. In the second place
mentioned above there is a misprint for the index of $\bar{P}$ which,
instead of $\max \left\{ i,j\right\} +1,$ should be $\max \left\{
i,j\right\} $.
\end{remark}

\section{Some New Inequalities}

The following result holds

\begin{theorem}
\label{t3.1} Let $\left( X,\left\Vert .\right\Vert \right) $ be a normed
linear space over the real or complex number field $\mathbb{K}$, $\mathbf{a}%
=\left( a_{1},...,a_{n}\right) \in \mathbb{K}^{n},\mathbf{p}=\left(
p_{1},...,p_{n}\right) \in \mathbb{R}^{n}$ and $\mathbf{x}=\left(
x_{1},...,x_{n}\right) \in X^{n}.$ Then one has the inequalities
\begin{eqnarray}
&&\left\Vert T_{n}\left( \mathbf{p};\mathbf{a},\mathbf{x}\right) \right\Vert
\label{3.1} \\
&\leq &\left\{ 
\begin{array}{l}
\max_{1\leq i\leq n-1}\left\vert \det \left( 
\begin{array}{cc}
P_{i} & P_{n} \\ 
A_{i}\left( \mathbf{p}\right) & A_{n}\left( \mathbf{p}\right)
\end{array}
\right) \right\vert \cdot \sum_{j=1}^{n-1}\left\Vert \Delta x_{j}\right\Vert
; \\ 
\\ 
\left( \sum_{i=1}^{n-1}\left\vert \det \left( 
\begin{array}{cc}
P_{i} & P_{n} \\ 
A_{i}\left( \mathbf{p}\right) & A_{n}\left( \mathbf{p}\right)
\end{array}
\right) \right\vert ^{q}\right) ^{1/q}\cdot \left(
\sum_{j=1}^{n-1}\left\Vert \Delta x_{j}\right\Vert ^{p}\right) ^{1/p}\text{ }
\\ 
\\ 
\text{for }p>1,\frac{1}{p}+\frac{1}{q}=1; \\ 
\\ 
\sum_{i=1}^{n-1}\left\vert \det \left( 
\begin{array}{cc}
P_{i} & P_{n} \\ 
A_{i}\left( \mathbf{p}\right) & A_{n}\left( \mathbf{p}\right)
\end{array}
\right) \right\vert \cdot \max_{1\leq j\leq n-1}\left\Vert \Delta
x_{j}\right\Vert .
\end{array}
\right.  \notag
\end{eqnarray}
All the inequalities in $\left( \ref{3.1}\right) $ are sharp in the sense
that the constants $1$ cannot be replaced by smaller constants.
\end{theorem}

\begin{proof}
Using the first identity in $\left( \ref{2.1.a}\right) $, we have
\begin{equation*}
\left\Vert T_{n}\left( \mathbf{p};\mathbf{a},\mathbf{x}\right) \right\Vert
\leq \sum_{i=1}^{n}\left\vert \det \left( 
\begin{array}{cc}
P_{i} & P_{n} \\ 
A_{i}\left( \mathbf{p}\right) & A_{n}\left( \mathbf{p}\right)
\end{array}
\right) \right\vert \left\Vert \Delta x_{i}\right\Vert .
\end{equation*}
Using H\"{o}lder's inequality, we deduce the desired result $\left( \ref{3.1}%
\right) .$

Let prove, for instance, that the constant $1$ in the second inequality is
best possible.

Assume, for $C>0$, we have that 
\begin{equation}
\left\Vert T_{n}\left( \mathbf{p};\mathbf{a},\mathbf{x}\right) \right\Vert
\leq C\left( \sum_{i=1}^{n-1}\left\vert \det \left( 
\begin{array}{cc}
P_{i} & P_{n} \\ 
A_{i}\left( \mathbf{p}\right) & A_{n}\left( \mathbf{p}\right)
\end{array}
\right) \right\vert ^{q}\right) ^{1/q}\left( \sum_{j=1}^{n-1}\left\Vert
\Delta x_{j}\right\Vert ^{p}\right) ^{1/p}  \label{3.2}
\end{equation}
for $p>1,\frac{1}{p}+\frac{1}{q}=1,n\geq 2.$

If we choose $n=2,$ then we get 
\begin{equation*}
T_{2}\left( \mathbf{p};\mathbf{a},\mathbf{x}\right) =p_{1}p_{2}\left(
a_{2}-a_{1}\right) \left( x_{2}-x_{1}\right) .
\end{equation*}
Also, for $n=2,$

\begin{equation*}
\left( \sum_{i=1}^{n-1}\left\vert \det \left( 
\begin{array}{cc}
P_{i} & P_{n} \\ 
A_{i}\left( \mathbf{p}\right) & A_{n}\left( \mathbf{p}\right)
\end{array}
\right) \right\vert ^{q}\right) ^{1/q}=\left\vert p_{1}p_{2}\right\vert
\left\vert a_{2}-a_{1}\right\vert
\end{equation*}
and
\begin{equation*}
\left( \sum_{j=1}^{n-1}\left\Vert \Delta x_{j}\right\Vert ^{p}\right)
^{1/p}=\left\Vert x_{2}-x_{1}\right\Vert .
\end{equation*}
Then by $\left( \ref{3.2}\right) $, holding for $n=2,p_{1},p_{2}>0,a_{1}\neq
a_{2},x_{2}\neq x_{1},$ we deduce $C\geq 1,$ proving that $1$ is the best
possible constant in that inequality.
\end{proof}

The following corollary for the uniform distribution of the probability $%
\mathbf{p}$ holds.

\begin{corollary}
\label{c3.2} With the assumptions of Theorem \ref{t3.1} for $\mathbf{a}$ and 
$\mathbf{x}$, we have the inequalities 
\begin{eqnarray*}
&&\left\Vert T_{n}\left( \mathbf{a},\mathbf{x}\right) \right\Vert \\
&\leq &\frac{1}{n^{2}}\times \left\{ 
\begin{array}{l}
\max_{1\leq i\leq n-1}\left\vert \det \left( 
\begin{array}{cc}
i & n \\ 
\sum_{k=1}^{i}a_{k} & \sum_{k=1}^{n}a_{k}
\end{array}
\right) \right\vert \cdot \sum_{j=1}^{n-1}\left\Vert \Delta x_{j}\right\Vert
; \\ 
\\ 
\left( \sum_{i=1}^{n-1}\left\vert \det \left( 
\begin{array}{cc}
i & n \\ 
\sum_{k=1}^{i}a_{k} & \sum_{k=1}^{n}a_{k}
\end{array}
\right) \right\vert ^{q}\right) ^{1/q}\cdot \left(
\sum_{j=1}^{n-1}\left\Vert \Delta x_{j}\right\Vert ^{p}\right) ^{1/p}\text{ }
\\ 
\\ 
\text{for }p>1,\frac{1}{p}+\frac{1}{q}=1; \\ 
\\ 
\sum_{i=1}^{n-1}\left\vert \det \left( 
\begin{array}{cc}
i & n \\ 
\sum_{k=1}^{i}a_{k} & \sum_{k=1}^{n}a_{k}
\end{array}
\right) \right\vert \cdot \max_{1\leq j\leq n-1}\left\Vert \Delta
x_{j}\right\Vert .
\end{array}
\right.
\end{eqnarray*}
\end{corollary}

The following result may be stated as well.

\begin{theorem}
\label{t3.3} With the assumptions of Theorem \ref{t3.1} and if $P_{i}\neq
0\left( i=1,...,n\right) ,$ then we have the inequalities
\begin{eqnarray}
&&\left\| T_{n}\left( \mathbf{p};\mathbf{a},\mathbf{x}\right) \right\| 
\label{3.4} \\
&\leq &\left| P_{n}\right| \times \left\{ 
\begin{array}{l}
\max_{1\leq i\leq n-1}\left| \frac{A_{n}\left( \mathbf{p}\right) }{P_{n}}-%
\frac{A_{i}\left( \mathbf{p}\right) }{Pi}\right| \cdot
\sum_{i=1}^{n-1}\left| P_{i}\right| \left\| \Delta x_{i}\right\| ; \\ 
\\ 
\left( \sum_{i=1}^{n-1}\left| P_{i}\right| \left| \frac{A_{n}\left( \mathbf{p%
}\right) }{P_{n}}-\frac{A_{i}\left( \mathbf{p}\right) }{Pi}\right|
^{q}\right) ^{1/q}\cdot \left( \sum_{i=1}^{n-1}\left| P_{i}\right| \left\|
\Delta x_{i}\right\| \right) ^{1/p} \\ 
\\ 
\text{for }p>1,\frac{1}{p}+\frac{1}{q}=1; \\ 
\\ 
\sum_{i=1}^{n-1}\left| P_{i}\right| \left| \frac{A_{n}\left( \mathbf{p}%
\right) }{P_{n}}-\frac{A_{i}\left( \mathbf{p}\right) }{Pi}\right| \cdot
\max_{1\leq i\leq n-1}\left\| \Delta x_{i}\right\| .
\end{array}
\right.   \notag
\end{eqnarray}
All the inequalities in $\left( \ref{3.4}\right) $ are sharp in the sense
that the constant $1$ cannot be replaced by a smaller constant.
\end{theorem}

\begin{proof}
Follows by the second identity in $\left( \ref{2.1.a}\right) $ and taking
into account that 
\begin{equation*}
\left\Vert T_{n}\left( \mathbf{p};\mathbf{a},\mathbf{x}\right) \right\Vert
\leq \left\vert P_{n}\right\vert \sum_{i=1}^{n-1}\left\vert \frac{%
A_{n}\left( \mathbf{p}\right) }{P_{n}}-\frac{A_{i}\left( \mathbf{p}\right) }{%
Pi}\right\vert \cdot \left\vert P_{i}\right\vert \left\Vert \Delta
x_{i}\right\Vert .
\end{equation*}
Using H\"{o}lder's weighted inequality, \ we easily deduce $\left( \ref{3.4}%
\right) .$

The sharpness of the constant may be shown in a similar manner. We omit the
details.
\end{proof}

The following corollary containing the unweighted inequalities holds.

\begin{corollary}
\label{c3.4} With the above assumptions for $\mathbf{a}$ and $\mathbf{x}$
one, has
\begin{eqnarray}
&&\left\| T_{n}\left( \mathbf{a},\mathbf{x}\right) \right\|   \label{3.5} \\
&\leq &\frac{1}{n}\times \left\{ 
\begin{array}{l}
\max_{1\leq i\leq n-1}\left| \frac{1}{n}\sum_{k=1}^{n}a_{k}-\frac{1}{i}%
\sum_{k=1}^{i}a_{k}\right| \cdot \sum_{i=1}^{n-1}i\left\| \Delta
x_{i}\right\| ; \\ 
\\ 
\left( \sum_{i=1}^{n-1}i\left| \frac{1}{n}\sum_{k=1}^{n}a_{k}-\frac{1}{i}%
\sum_{k=1}^{i}a_{k}\right| ^{q}\right) ^{1/q}\cdot \left(
\sum_{i=1}^{n-1}i\left\| \Delta x_{i}\right\| ^{p}\right) ^{1/p} \\ 
\\ 
\text{for }p>1,\frac{1}{p}+\frac{1}{q}=1; \\ 
\\ 
\sum_{i=1}^{n-1}i\left| \frac{1}{n}\sum_{k=1}^{n}a_{k}-\frac{1}{i}%
\sum_{k=1}^{i}a_{k}\right| \cdot \max_{1\leq i\leq n-1}\left\| \Delta
x_{i}\right\| .
\end{array}
\right.   \notag
\end{eqnarray}
The inequalities in $\left( \ref{3.5}\right) $ are sharp in the sense
mentioned above.
\end{corollary}

Another type of inequalities may be stated if one uses the third identity in 
$\left( \ref{2.1.a}\right) $.

\begin{theorem}
\label{t3.4} With the assumptions in Theorem \ref{t3.1} and if $P_{i},%
\overline{P}_{i}\neq 0,$ $i\in \left\{ 1,...,n-1\right\} ,$ then we have the
inequalities
\begin{eqnarray}
&&\left\Vert T_{n}\left( \mathbf{p};\mathbf{a},\mathbf{x}\right) \right\Vert
\label{3.6} \\
&\leq &\left\{ 
\begin{array}{l}
\max_{1\leq i\leq n-1}\left\vert \frac{\overline{A}_{i}\left( \mathbf{p}%
\right) }{\overline{P}_{i}}-\frac{A_{i}\left( \mathbf{p}\right) }{Pi}%
\right\vert \cdot \sum_{i=1}^{n-1}\left\vert P_{i}\right\vert \left\vert 
\overline{P}_{i}\right\vert \left\Vert \Delta x_{i}\right\Vert ; \\ 
\\ 
\left( \sum_{i=1}^{n-1}\left\vert P_{i}\right\vert \left\vert \overline{P}%
_{i}\right\vert \left\vert \frac{\overline{A}_{i}\left( \mathbf{p}\right) }{%
\overline{P}_{i}}-\frac{A_{i}\left( \mathbf{p}\right) }{Pi}\right\vert
^{q}\right) ^{1/q}\cdot \left( \sum_{i=1}^{n-1}\left\vert P_{i}\right\vert
\left\vert \overline{P}_{i}\right\vert \left\Vert \Delta x_{i}\right\Vert
^{p}\right) ^{1/p} \\ 
\\ 
\text{for }p>1,\frac{1}{p}+\frac{1}{q}=1; \\ 
\\ 
\sum_{i=1}^{n-1}\left\vert P_{i}\right\vert \left\vert \overline{P}%
_{i}\right\vert \left\vert \frac{\overline{A}_{i}\left( \mathbf{p}\right) }{%
\overline{P}_{i}}-\frac{A_{i}\left( \mathbf{p}\right) }{Pi}\right\vert \cdot
\max_{1\leq i\leq n-1}\left\Vert \Delta x_{i}\right\Vert .
\end{array}
\right.  \notag
\end{eqnarray}
In particular, if $p_{i}=\frac{1}{n},i\in \left\{ 1,...,n\right\} ,$ then we
have 
\begin{eqnarray}
&&\left\Vert T_{n}\left( \mathbf{a},\mathbf{x}\right) \right\Vert
\label{3.7} \\
&\leq &\frac{1}{n^{2}}\cdot \left\{ 
\begin{array}{l}
\max_{1\leq i\leq n-1}\left\vert \frac{1}{n-i}\sum_{k=i+1}^{n}a_{k}-\frac{1}{%
i}\sum_{k=1}^{i}a_{k}\right\vert \cdot \sum_{i=1}^{n-1}i\left( n-i\right)
\left\Vert \Delta x_{i}\right\Vert ; \\ 
\\ 
\left( \sum_{i=1}^{n-1}i\left( n-i\right) \left\vert \frac{1}{n-i}%
\sum_{k=i+1}^{n}a_{k}-\frac{1}{i}\sum_{k=1}^{i}a_{k}\right\vert ^{q}\right)
^{1/q} \\ 
\\ 
\times \left( \sum_{i=1}^{n-1}i\left( n-i\right) \left\Vert \Delta
x_{i}\right\Vert ^{p}\right) ^{1/p}\text{ for }p>1,\frac{1}{p}+\frac{1}{q}=1;
\\ 
\\ 
\sum_{i=1}^{n-1}i\left( n-i\right) \left\vert \frac{1}{n-i}%
\sum_{k=i+1}^{n}a_{k}-\frac{1}{i}\sum_{k=1}^{i}a_{k}\right\vert \cdot
\max_{1\leq i\leq n-1}\left\Vert \Delta x_{i}\right\Vert .
\end{array}
\right.  \notag
\end{eqnarray}
The inequalities in $\left( \ref{3.6}\right) $and $\left( \ref{3.7}\right) $
are sharp in the above mentioned sense.
\end{theorem}

A different approach may be considered if one uses the representation in
terms of double sums for the \v{C}eby\v{s}ev functional provided by the
Theorem \ref{t2.2}.

The following result holds.

\begin{theorem}
\label{t3.5}With the assumptions in Theorem \ref{t3.1}, we have the
inequalities
\begin{eqnarray}
&&\left\| T_{n}\left( \mathbf{p};\mathbf{a},\mathbf{x}\right) \right\| 
\label{3.8} \\
&\leq &\left\{ 
\begin{array}{l}
\max_{1\leq i,j\leq n-1}\left\{ \left| P_{\min \left\{ i,j\right\} }\right|
\left| \bar{P}_{\max \left\{ i,j\right\} }\right| \right\} \cdot
\sum_{i=1}^{n-1}\left| \Delta a_{i}\right| \sum_{i=1}^{n-1}\left\| \Delta
x_{i}\right\| ; \\ 
\\ 
\left( \sum_{i=1}^{n-1}\sum_{j=1}^{n-1}\left| P_{\min \left\{ i,j\right\}
}\right| ^{q}\left| \bar{P}_{\max \left\{ i,j\right\} }\right| ^{q}\right)
^{1/q} \\ 
\\ 
\times \left( \sum_{i=1}^{n-1}\left| \Delta a_{i}\right| ^{p}\right)
^{1/p}\left( \sum_{i=1}^{n-1}\left\| \Delta x_{i}\right\| ^{p}\right) ^{1/p}%
\text{ for }p>1,\frac{1}{p}+\frac{1}{q}=1; \\ 
\\ 
\sum_{i=1}^{n-1}\sum_{j=1}^{n-1}\left| P_{\min \left\{ i,j\right\} }\right|
\left| \bar{P}_{\max \left\{ i,j\right\} }\right|  \\ 
\\ 
\times \max_{1\leq i\leq n-1}\left| \Delta a_{i}\right| \max_{1\leq i\leq
n-1}\left\| \Delta x_{i}\right\| .
\end{array}
\right.   \notag
\end{eqnarray}
The inequalities are sharp in the sense mentioned above.
\end{theorem}

The proof follows by the identity $\left( \ref{2.7}\right) $ on using
H\"{o}lder's inequality for double sums and we omit the details.

Now, define
\begin{equation*}
k_{\infty }:=\max_{1\leq i,j\leq n-1}\left\{ \frac{\min \left\{ i,j\right\} 
}{n}\left( 1-\frac{\max \left\{ i,j\right\} }{n}\right) \right\} ,n\geq 2.
\end{equation*}
Using the elementary inequality
\begin{equation*}
ab\leq \frac{1}{4}\left( a+b\right) ^{2},\text{ \ }a,b\in R;
\end{equation*}
we deduce
\begin{eqnarray*}
\min \left\{ i,j\right\} \cdot \left( n-\max \left\{ i,j\right\} \right) 
&\leq &\frac{1}{4}\left( n+\min \left\{ i,j\right\} -\max \left\{
i,j\right\} \right) ^{2} \\
&=&\frac{1}{4}\left( n-\left| i-j\right| \right) ^{2},1\leq i,j\leq n-1.
\end{eqnarray*}
Consequently, we observe that
\begin{equation*}
k_{\infty }\leq \frac{1}{4n^{2}}\max_{1\leq i,j\leq n-1}\left\{ \left(
n-\left| i-j\right| \right) ^{2}\right\} =\frac{1}{4}.
\end{equation*}
We may state now the following corollary of Theorem \ref{t3.5}.

\begin{corollary}
\label{c3.6} Let $\left( X,\left\Vert .\right\Vert \right) $ be a normed
linear space, $\mathbf{a}=\left( a_{1},...,a_{n}\right) \in \mathbb{K}^{n}$
and $\mathbf{x}=\left( x_{1},...,x_{n}\right) \in \mathbb{X}^{n}.$ Then we
have the inequality
\begin{equation}
\left\Vert T_{n}\left( \mathbf{a},\mathbf{x}\right) \right\Vert \leq
k_{\infty }\sum_{i=1}^{n-1}\left\vert \Delta a_{i}\right\vert
\sum_{i=1}^{n-1}\left\Vert \Delta x_{i}\right\Vert \leq \frac{1}{4}%
\sum_{i=1}^{n-1}\left\vert \Delta a_{i}\right\vert
\sum_{i=1}^{n-1}\left\Vert \Delta x_{i}\right\Vert .  \label{3.9}
\end{equation}
The constant $\frac{1}{4}$ cannot be replaced in general by a smaller
constant.
\end{corollary}

\begin{remark}
\label{r3.7} The inequality $\left( \ref{3.9}\right) $ is better than the
second inequality in Corollary \ref{c1.2}.
\end{remark}

Consider now, for $q>1,$ the number
\begin{equation*}
k_{q}:=\frac{1}{n^{2}}\left( \sum_{i=1}^{n-1}\sum_{j=1}^{n-1}\left[ \min
\left\{ i,j\right\} \cdot \left( n-\max \left\{ i,j\right\} \right) \right]
^{q}\right) ^{1/q}.
\end{equation*}
We observe, by the symmetry of the terms under the sums symbol, we have that
\begin{equation*}
k_{q}=\frac{1}{n^{2}}\left( 2\sum_{1\leq i<j\leq n-1}i^{q}\left( n-j\right)
^{q}+\sum_{i=1}^{n-1}i^{q}\left( n-i\right) ^{q}\right) ^{1/q},
\end{equation*}
that may be computed exactly if $q=2$ or another natural number.

Since, as above,
\begin{equation*}
\left[ \min \left\{ i,j\right\} \cdot \left( n-\max \left\{ i,j\right\}
\right) \right] ^{q}\leq \frac{1}{4^{q}}\left( n-\left| i-j\right| \right)
^{2q}
\end{equation*}
we deduce
\begin{eqnarray*}
k_{q} &\leq &\frac{1}{4n^{2}}\left( \sum_{i=1}^{n-1}\sum_{j=1}^{n-1}\left(
n-\left| i-j\right| \right) ^{2q}\right) ^{1/q} \\
&\leq &\frac{1}{4n^{2}}\left[ \left( n-1\right) ^{2}n^{2q}\right] ^{1/q}=%
\frac{1}{4}\left( n-1\right) ^{2/q}.
\end{eqnarray*}
Consequently, we may state the following corollary as well.

\begin{corollary}
\label{c3.8} With the assumption in Corollary \ref{c3.6}, we have the
inequalities
\begin{eqnarray*}
\left\Vert T_{n}\left( \mathbf{a},\mathbf{x}\right) \right\Vert &\leq
&k_{q}\left( \sum_{i=1}^{n-1}\left\vert \Delta a_{i}\right\vert ^{p}\right)
^{1/p}\left( \sum_{i=1}^{n-1}\left\Vert \Delta x_{i}\right\Vert ^{p}\right)
^{1/p} \\
&\leq &\frac{1}{4}\left( n-1\right) ^{2/q}\left( \sum_{i=1}^{n-1}\left\vert
\Delta a_{i}\right\vert ^{p}\right) ^{1/p}\left( \sum_{i=1}^{n-1}\left\Vert
\Delta x_{i}\right\Vert ^{p}\right) ^{1/p};
\end{eqnarray*}
provided $p>1,\frac{1}{p}+\frac{1}{q}=1.$ The constant $\frac{1}{4}$ cannot
be replaced in general by a smaller constant.
\end{corollary}

Finally, if we denote
\begin{equation*}
k_{1}:=\frac{1}{n^{2}}\sum_{i=1}^{n-1}\sum_{j=1}^{n-1}\left[ \min \left\{
i,j\right\} \cdot \left( n-\max \left\{ i,j\right\} \right) \right] ,
\end{equation*}
then we observe, for $\mathbf{u}=\left( \frac{1}{n},...,\frac{1}{n}\right) ,%
\mathbf{e}=\left( 1,2,...,n\right) ,$ that 
\begin{equation*}
k_{1}=T_{n}\left( \mathbf{u};\mathbf{e},\mathbf{e}\right) =\frac{1}{n}%
\sum_{i=1}^{n}i^{2}-\left( \frac{1}{n}\sum_{i=1}^{n}i\right) ^{2}=\frac{1}{12%
}\left( n^{2}-1\right) ,
\end{equation*}
and by Theorem \ref{t3.5}, we deduce the inequality
\begin{equation*}
\left\| T_{n}\left( \mathbf{a},\mathbf{x}\right) \right\| \leq \frac{1}{12}%
\left( n^{2}-1\right) \max_{1\leq j\leq n-1}\left| \Delta a_{j}\right|
\max_{1\leq j\leq n-1}\left\| \Delta x_{j}\right\| .\text{ }
\end{equation*}
Note that, the above inequality, has been discovered with a different method
in \cite{DB}. The constant $\frac{1}{12}$ is best possible.

\end{document}